\documentclass[leqno]{article}
\usepackage[utf8]{inputenc}
\usepackage{systeme}
\usepackage{fullpage}
\usepackage[dvipsnames]{xcolor}
\usepackage{authblk}
\usepackage{setspace}

\usepackage{mathtools,amsfonts,amssymb,amsmath,amsthm}
\usepackage{algcompatible}
\usepackage{mathtools}
\usepackage{amssymb}
\usepackage{cases}

\newtheorem{theorem}{Theorem}
\newtheorem{prop}{Proposition}

\newtheorem{lemma}{Lemma}
\newtheorem{remark}{Remark}
\newcommand{\FF}{\mathbb{F}}

\newcommand{\N}{\mathbb{N}}

\newcommand{\R}{\mathbb{R}}

\usepackage{pgf,tikz,pgfplots}
\usepackage{tikz-3dplot}
\usetikzlibrary{shapes,arrows.meta,bending,positioning,shadows,patterns,pgfplots.groupplots,decorations.markings,chains}
\usepackage{pgfplots}
\usetikzlibrary{calc}
\usetikzlibrary{shapes.geometric}
\usetikzlibrary{arrows.meta}
\pgfplotsset{compat=newest}
\usepackage{mathrsfs}

\usepackage{color,ulem}
\usepackage{graphicx}%
\usepackage[active]{srcltx}
\usepackage{mathabx}
\usepackage{verbatim}
\usepackage{cancel}
\usepackage{tikz,pgfplots}
\usepackage{tabu}

\usepackage{hyperref}

\usepackage{caption}
\usepackage{subcaption}

\usepackage[official]{eurosym}

\usepackage[boxruled]{algorithm2e} 
\SetKwRepeat{Do}{do}{while}

\def \R{\mathbb R}

\def \N{\mathbb N}

\newcommand\p\varphi
\newcommand\e\epsilon

\title{An elementary approach based on variational inequalities for modelling a friction-based locomotion problem}

\author[1]{Panyu Chen}
\author[2,3]{Álvaro Mateos González}
\author[4]{Laurent Mertz}
\affil[1]{Department of Computer Science, Duke University}
\affil[2]{University of Michigan-Shanghai Jiao Tong University Joint Institute, Shanghai Jiao Tong University, Shanghai, 200240, China}
\affil[3]{Modeling, Data Analysis \& Computational Tools for Biology Research Group, Biomathematics Unit, Department of Biodiversity, Ecology \& Evolution, Faculty of Biological Sciences, Complutense University of Madrid, 28040 Madrid, Spain}
\affil[4]{Department of Mathematics, City University of Hong Kong, Kowloon, Hong Kong, China}

\begin{document}

\maketitle
\begin{abstract}
We propose an elementary proof based on a penalization technique to show the existence and uniqueness of the solution to a system of variational inequalities modelling the friction-based motion of a two-body crawling system. Here for each body, the static and dynamic friction coefficients are equal.
\end{abstract}

\section{Introduction}
\label{sec:scallop}
The ability of a living organism to move from one place to another is essential to its survival. Indeed, it makes it possible to meet many needs such as finding food, a habitat and a mate, or escaping predators and danger. Living organisms use a broad class of locomotion strategies such as swimming in the water, flying through the air or running, crawling, and slithering using friction between their body and the ground on land \cite{alexander2003principles}. The understanding and simulation of the mechanisms underlying these locomotion strategies have promising applications in robotics \cite{Figurina04,hirose1993biologically,transeth_pettersen_liljebäck_2009} by biomimetics. In the present work, we propose a system of variational inequalities to model friction-based locomotion for multibody systems using friction forces between their bodies and a surface.

\subsection{Physical formulation of a two-body crawling system using dry friction}
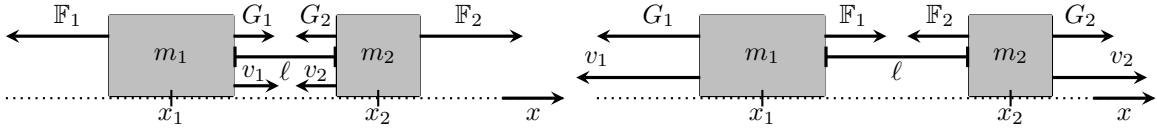
\begin{figure}[h!]
\begin{center} 
\begin{tikzpicture}[scale=0.55]

\coordinate (GL) at (-4,0); 
\coordinate (GR) at (9.5,0); 

\coordinate (1BL) at (-1.5,0.05); 
\coordinate (1BR) at (1.5, 0.05); 
\coordinate (1TL) at (-1.5, 2); 
\coordinate (1TR) at (1.5, 2); 
\coordinate (1C) at (0,1); 

\coordinate (2BL) at (4,0.05); 
\coordinate (2BR) at (6, 0.05); 
\coordinate (2TL) at (4, 2); 
\coordinate (2TR) at (6, 2); 
\coordinate (2C) at (5, 1); 

\coordinate (forceF1) at (-2.5,2.0); 
\coordinate (forceF2) at (7.2,2.0); 
\coordinate (forceG1) at (2.1,2.0); 
\coordinate (forceG2) at (3.5,2.0); 
\coordinate (movement) at (8.75,-0.4); 
\coordinate (velocityV1) at (2,0.6); 
\coordinate (velocityV2) at (3.5,0.6); 

\draw[dotted,thick] (GL) -- (GR);			
\draw[thick] (-1.5,0.05) -- (-1.5,2);		
\draw[thick] (1.5,0.05) -- (1.5,2);			
\draw[thick] (-1.5,0.05) -- (1.5,0.05);		
\draw[thick] (-1.5,2) -- (1.5,2);			
\draw[thick] (4.0,0.05) -- (4.0,2);			
\draw[thick] (6.0,0.05) -- (6.0,2);			
\draw[thick] (4.0,0.05) -- (6.0,0.05);		
\draw[thick] (4.0,2) -- (6.0,2);			

\fill[color=gray!50!white] (1BL) rectangle (1TR);		
\fill[color=gray!50!white] (2BL) rectangle (2TR);		

\draw[|-|,>=stealth,very thick,black] (1.5,1.0) -> (4.0,1.0);		
\draw[->,>=stealth,very thick,black] (-1.5,1.5) -> (-4.0,1.5);	
\draw[->,>=stealth,very thick,black] (4.0,1.5) -> (3,1.5);	
\draw[->,>=stealth,very thick,black] (6,1.5) -> (8.5,1.5);	
\draw[->,>=stealth,very thick,black] (1.5,1.5) -> (2.5,1.5);	
\draw[->,>=stealth,very thick,black] (1.5,0.3) -> (2.6,0.3);	
\draw[->,>=stealth,very thick,black] (4.0,0.3) -> (3.0,0.3);	
\draw[->,>=stealth,very thick,black] (8,0) -> (9.5,0);	

\draw[-,>=stealth,thick,black] (0,-0.2) -> (0,0.2);	
\draw[-,>=stealth,thick,black] (5,-0.2) -> (5,0.2);	

\node at (forceF1) {$\FF_1$};
\node at (forceF2) {$\FF_2$};
\node at (forceG1) {$G_1$};
\node at (forceG2) {$G_2$};
\node at (velocityV1) {$v_1$};
\node at (velocityV2) {$v_2$};
\node at (movement) {$x$};

 \node at (barycentric cs:1BL=1,1TL=1,1TR=1,1BR=1) {$m_1$};
 \node at (barycentric cs:1BL=1,1TL=-0.2,1TR=-0.2,1BR=1) {$x_1$};
 \node at (barycentric cs:2BL=1,2TL=1,2TR=1,2BR=1) {$m_2$};
 \node at (barycentric cs:2BL=1,2TL=-0.2,2TR=-0.2,2BR=1) {$x_2$};

 \node at (2.75,0.65) {$\ell$};
 
\end{tikzpicture}
\begin{tikzpicture}[scale=0.55]

\coordinate (GL) at (-4,0); 
\coordinate (GR) at (9.5,0); 

\coordinate (1BL) at (-1.5,0.05); 
\coordinate (1BR) at (1.5, 0.05); 
\coordinate (1TL) at (-1.5, 2); 
\coordinate (1TR) at (1.5, 2); 
\coordinate (1C) at (0,1); 

\coordinate (2BL) at (5,0.05); 
\coordinate (2BR) at (7, 0.05); 
\coordinate (2TL) at (5, 2); 
\coordinate (2TR) at (7, 2); 
\coordinate (2C) at (6, 1); 

\coordinate (forceG1) at (-2.5,2.0); 
\coordinate (forceG2) at (7.7,2.0); 
\coordinate (forceF1) at (2.2,2.0); 
\coordinate (forceF2) at (4.3,2.0); 
\coordinate (movement) at (8.75,-0.4); 
\coordinate (velocityV1) at (-4,0.95); 
\coordinate (velocityV2) at (8.7,0.95); 

\draw[dotted,thick] (GL) -- (GR);			
\draw[thick] (-1.5,0.05) -- (-1.5,2);		
\draw[thick] (1.5,0.05) -- (1.5,2);			
\draw[thick] (-1.5,0.05) -- (1.5,0.05);		
\draw[thick] (-1.5,2) -- (1.5,2);			
\draw[thick] (5.0,0.05) -- (5.0,2);			
\draw[thick] (7.0,0.05) -- (7.0,2);			
\draw[thick] (5.0,0.05) -- (7.0,0.05);		
\draw[thick] (5.0,2) -- (7.0,2);			

\fill[color=gray!50!white] (1BL) rectangle (1TR);		
\fill[color=gray!50!white] (2BL) rectangle (2TR);		

\draw[|-|,>=stealth,very thick,black] (1.5,1.0) -> (5.0,1.0);		
\draw[->,>=stealth,very thick,black] (-1.5,1.5) -> (-4.0,1.5);	
\draw[->,>=stealth,very thick,black] (5.0,1.5) -> (3.5,1.5);	
\draw[->,>=stealth,very thick,black] (7,1.5) -> (8.5,1.5);	
\draw[->,>=stealth,very thick,black] (1.5,1.5) -> (3.0,1.5);	
\draw[->,>=stealth,very thick,black] (-1.5,0.5) -> (-4.5,0.5);	
\draw[->,>=stealth,very thick,black] (7,0.5) -> (9.3,0.5);	
\draw[->,>=stealth,very thick,black] (8,0) -> (9.5,0);	

\draw[-,>=stealth,thick,black] (0,-0.2) -> (0,0.2);	
\draw[-,>=stealth,thick,black] (6,-0.2) -> (6,0.2);	

\node at (forceF1) {$\FF_1$};
\node at (forceF2) {$\FF_2$};
\node at (forceG1) {$G_1$};
\node at (forceG2) {$G_2$};
\node at (movement) {$x$};
\node at (velocityV1) {$v_1$};
\node at (velocityV2) {$v_2$};

\node at (barycentric cs:1BL=1,1TL=1,1TR=1,1BR=1) {$m_1$};
\node at (barycentric cs:1BL=1,1TL=-0.2,1TR=-0.2,1BR=1) {$x_1$};
\node at (barycentric cs:2BL=1,2TL=1,2TR=1,2BR=1) {$m_2$};
\node at (barycentric cs:2BL=1,2TL=-0.2,2TR=-0.2,2BR=1) {$x_2$};

\node at (3.25,0.65) {$\ell$};
 
\end{tikzpicture}
\hspace{1cm}
\end{center}

\caption{Bodies of masses $m_1$ and $m_2$ at respective positions $x_1$ and $x_2$. The linkage $\ell = x_2 - x_1$ imposes contact forces $G_1 = - G_2$ on the bodies. The directions of friction forces $\FF_1$ and $\FF_2$ oppose the body movements. 
Left: contraction phase. Right: extension phase.    
\label{fig:friction}}
\end{figure}

We consider a setting similar to that of \cite{wagnerlauga13}: a system composed of two bodies with masses $m_1$ and $m_2$ lying on a flat surface. We assume all motion of the system is uni-dimensional along a single axis, thus we can represent the positions of masses $m_1$ and $m_2$ by two real numbers $x_1$ and $x_2$.
The velocity and the acceleration of mass $m_i$, $i \in \{1, 2\}$, are respectively denoted by $\dot x_i = \frac{{\rm d} x_i}{{\rm d}t}$ and $\ddot x_i$. 
A mechanism (linkage) connects the two bodies and imposes the temporal variation of their distance $\ell(t) \triangleq x_2(t) - x_1(t)$. 
For instance, if $\ell$ is a given oscillating function, it corresponds to an alternation between extension phases in which the bodies are pushed apart and contraction phases in which they are pulled closer. In this way, an overall motion of the system can be generated by the 
dry friction forces. 
The motion of the two-body system is described by Newton's second law: for each $i = 1,2$, we have 
$m_i \ddot x_i + \mathbb{F}_i = G_i$, where $G_i$ is a contact force (push or pull) applied to the body of mass $m_i$ which is induced by $\ell$ and the body with mass $m_j, j \neq i$ and $\mathbb{F}_i$ is a dry friction force described by Coulomb's model.
The two forces $\mathbb{F}_i$ and $G_i$ are coupled and described as follows:
\[
 \mathbb{F}_i \triangleq 
\left\{
 \begin{array}{ll}
 G_i, \: & \mbox{when} \: \dot x_i =0 \mbox{ and } |G_i| \leq f_i, \medskip\\
 \sigma_i f_i, \: & \mbox{when} \: \dot x_i \neq0 \; \mbox{or} \: (\dot x_i = 0 \mbox{ and } |G_i| > f_i)
\end{array}
\right.
\]
where $\sigma_i \triangleq \textup{sign}(\dot x_i)$ when $\dot x_i \neq 0$ or $\sigma_i \triangleq \textup{sign}(G_i)$ when $\dot x_i = 0$. 
Here $f_i$ is the dry friction coefficient associated with the mass $m_i$. 
See Figure~\ref{fig:friction} for illustration.
Newton's third law implies that 
$G_1 = -G_2$,
and then using $x_2 = x_1 + \ell$ to reduce the aforementioned equations of motion we obtain
$$
G_2 = \frac{m_1m_2}{M} \left ( \ddot \ell + \frac{\mathbb{F}_2}{m_2} - \frac{\mathbb{F}_1}{m_1} \right )
$$
where $M \triangleq m_1+m_2$.
As a consequence, the equations of motion can be formulated as 
\begin{equation} \label{eq:Newton2}
M \ddot x_1 + \mathbb{F}_1 = - m_2 \ddot \ell - \mathbb{F}_2, \qquad 
M \ddot x_2 + \mathbb{F}_2 = m_1 \ddot \ell - \mathbb{F}_1.  
\end{equation}
\noindent
Note that the equation of $x_2$ is redundant since $x_2$ is entirely determined by $x_2 = x_1 + \ell$, where $\ell$ is a known input. Thus, we are concerned with clarifying the dynamics of $x_1$ in order to represent the whole system.

\subsection{Review of methods for friction-based locomotion}

The dry friction forces in our paper describing solid-solid friction are rather standard, see chapter 10 of \cite{Popov2010}. In a simple setting, the response of a single mass subjected to dry friction and external forces with appropriate regularity (possibly random) is well understood in terms of standard variational inequalities or differential inclusions \cite{Brezis73,moreau1977application,pardoux2014stochastic}.  
Here, our setting with two masses is inspired by the work \cite{wagnerlauga13}, where the authors prove that locomotion can occur in the framework introduced above provided the dry friction coefficients $f_1$ and $f_2$ are distinct. The time-variation of the body–body separation varies periodically, but asymmetrically in time. The authors study stick and slip phases separately and provide numerical simulations and asymptotic interpretations of their results.  
Models related to ours have recurrently been used in robotics inspired by biology, particularly for modular snake robots able to navigate complex environments. This ability spurred a theoretical and practical interest in modeling, building and potentially using such robots in impractical terrain. Starting in the 70s, Hirose and coauthors designed models for the biomechanics of snake locomotion, focusing on the coordination of contraction/extension of a chain of modular segments and the management of dry friction forces necessary for effective locomotion on rough surfaces \cite{hirose1993biologically,UmetaniHirose1974}. A survey on snake robot models can be found in \cite{transeth_pettersen_liljebäck_2009}. More recently, several research groups have built experimental prototypes to study settings very closely related to ours. In \cite{behn2017dynamics}, the authors model the 1-dimensional motion of a system of three connected masses on the dry frictional plane, enabled by periodic changes in inter-mass distances. Those distances are explicitly defined as continuous piece-wise parabolic periodic functions. The motion equations are solved for those analytic expressions and implemented in a prototype. In \cite{bolotnik2018periodic}, the authors model a two-body crawler on an inclined plane under dry friction. Again, the system is actuated by periodic changes in distance masses. They analyze numerically and empirically the range of inclination angles and excitation laws that enable the system to crawl upwards. Beyond the case of unidimensional motion, several works have been focused on modeling 2- or 3-dimensional motion. In \cite{Figurina04}, the quasi-static motion of a two-link system on a frictional 2d plane is investigated. In \cite{marvi2011scalybot}, the development of a snake-inspired robot with active control of friction showcases practical applications of such models, enabling efficient movement across diverse terrains. On the theoretical analysis side, Gidoni and coauthors have studied abstract rate-independent systems motivated by applications to soft and elastic crawlers \cite{gidoni2018rate,Gidoni17,Gidoni2021}. Their work in \cite{Gidoni2023} is more general than ours, but their setup is more complex.  
In our work, we remain in the setup of \cite{wagnerlauga13}, allowing for an elementary approach based on a penalization technique to obtain the existence uniqueness of a solution to Equation \eqref{eq:Newton2}.

\subsection{Our contribution: a system of variational inequalities modelling a two-body crawling system}
For $i = 1,2$, we define the functions $\varphi_i(y) = f_i |y|$, for $y \in \R$. Given $\ell(.)$ twice continuously differentiable on $[0,\infty)$ and $y_0 \in \mathbb{R}$, we consider the mathematical problem of finding three functions $y,k_1,k_2: [0,\infty) \to \mathbb{R}$ such that 
\begin{align*}
\label{myproblem}
\tag{$\mathcal{S}_2$}
(i) \qquad & y \in \mathcal{C}([0,\infty);\mathbb{R}), y(0) = y_0,\\
(ii) \qquad & k_1,k_2 \in \mathcal{C}([0,\infty);\mathbb{R}) \cap \mathcal{BV}_{loc}([0,\infty);\mathbb{R}), \: {\text{with}} \: k_1(0) = 0,\\
(iii) \qquad & \forall t \geq 0, \: M y(t) + k_1(t) + k_2(t) = - m_2 \dot \ell(t),\\
(iv) \qquad &  \forall t \geq 0,
\: \forall \: 0 \leq s \leq t, 
\: \forall \mathfrak{z} \in \mathcal{C}([s,t];\mathbb{R}),\\
& \int_s^t [\mathfrak{z}(r) - y(r)] \textup{d} k_1(r) + \int_s^t \varphi_1(y(r)) \textup{d} r \leq \int_s^t \varphi_1 (\mathfrak{z}(r)) \textup{d} r,\\
(v) \qquad & \forall t \geq 0,
\: \forall \: 0 \leq s \leq t, 
\: \forall \mathfrak{z} \in \mathcal{C}([s,t];\mathbb{R}),\\
& \int_{s}^{t} [\mathfrak{z}(r) - y(r) - \dot \ell (r)] \textup{d} k_2(r) + \int_{s}^{t} \varphi_2 \left( y(r) + \dot \ell(r) \right) \textup{d} r \leq \int_s^t \varphi_2 (\mathfrak{z}(r)) \textup{d} r.
\end{align*}
Due to the presence of a variational parameter $\mathfrak z$, the fourth and fifth conditions are also called variational inequalities (VI).  Using the notation of \cite{pardoux2014stochastic}, 
the first VI above corresponds to the differential inclusion $\textup{d} k_1(t) \in \partial \varphi_1(y(t))(\textup{d} t)$ whereas the second one correspond to $\textup{d} k_2(t) \in \partial \varphi_2(y(t) + \dot \ell(t))(\textup{d} t)$. Here we have $\partial \varphi_i(y) = \{ f_i \textup{sign}(y) \} \: \mbox{ if } \: y \neq 0$ or $[-f_i,f_i]$ if $y=0$.
As a consequence, the problem above clarifies the sense of the following nonstandard multivalued equation
$M \textup{d} y + \partial \varphi_1(y)(\textup{d} t) + \partial \varphi_2(y+\dot \ell)(\textup{d} t) \ni - m_2 \dot \ell \textup{d} t$.
Then the definition of $x_1$ becomes straightforward
\[
\forall t \geq 0, \: x_1(t) = x_{1,0} + \int_0^t y(r) \textup{d} r.
\]
\begin{theorem} \label{thm:EUsol}
The problem \eqref{myproblem} has a unique solution. 
\end{theorem}
\begin{remark}
Our first theorem above resolves a nonstandard problem of Skorokhod. If we remove $\partial \varphi_1(y+\ell)(\textup{d} t)$ in the multivalued equation above then it corresponds to a standard Skorokhod problem as shown on page 231 of \cite{pardoux2014stochastic} for instance. In terms of friction forces, roughly speaking, $\partial \varphi_1(y)$ and $\partial \varphi_2(y+\dot \ell)$ corresponds to $\mathbb{F}_1$ and $\mathbb{F}_2$ respectively.
It is worth mentioning that the same physical problem with $p \geq 2$ aligned masses can be treated in our mathematical framework with a straightforward generalization of the case $p=2$.
\end{remark}

\section{Proof of Theorem~\ref{thm:EUsol}} 
The proof is composed of four steps. The first step consists in studying a regularised version of \eqref{myproblem} using the Moreau-Yosida regularization with two parameters $n_1$ and $n_2$. The regularized equation whose state variable is $y^n$ with $n=(n_1,n_2)$ can be analyzed using standard method of ODEs. In the second step, we show the sequence $y^n$ has the Cauchy property with respect to the supremum norm $\| y \| := \sup_{t \in \R_+} |y(t)|$. As a consequence, a limit exists.
Then, in the third step, we identify the limit as a solution to a system of variational inequalities. Finally, in the last step, we prove the uniqueness of the limit.  

\subsection{Moreau-Yosida regularisation}
The Moreau-Yosida regularisation of $\varphi_i$, which is a standard tool in non-smooth convex analysis in~\cite[Chap.2]{Brezis73}, is defined as follows: for each $n_i \in \mathbb{N}^\star$, for each $y \in \mathbb{R}$,
$$
\varphi_{i,n_i}(y) \triangleq \inf \limits_{y' \in \mathbb{R}} \left \{ \frac{n_i}{2} (y-y')^2 + \varphi_i(y') \right \}
= 
\begin{cases}
f_i \frac{n_i}{2} y^2, & \text{ if } n_i |y| \leq f_i,\\
f_i \left ( |y| - \frac{1}{2n_i} \right ), & \text{ if } n_i |y| > f_i. 
\end{cases}
$$
Thus, $\lim \limits_{n_i \to \infty} \varphi_{i,n_i}(y) = \varphi_i(y)$. 
For each $n_i \in \mathbb{N}^\star$, $\varphi_{i,n_i}(y)$ is differentiable and $\varphi_{i,n_i}'(y) = {\rm proj}_{[-f_i,f_i]}(n_i y)$.
In this context, a regularisation of the problem \eqref{myproblem} consists in replacing $\varphi_i$ by $\varphi_{i,n_i}$. 
Here given $\ell(.)$ twice continuously differentiable on $[0,\infty)$ and $y_0 \in \mathbb{R}$,  we want to find a function $y^n \in C([0,\infty);\mathbb{R})$ where $n \triangleq (n_1,n_2)$ such that 
\begin{equation}
\label{eq:smoothintegrated}
    \forall t \geq 0, \: M y^n(t) 
+ \int_0^t \varphi_{1,n_1}'(y^n(s)) \textup{d} s 
+ \int_0^t \varphi_{2,n_2}'(y^n(s) + \dot \ell(s)) \textup{d} s 
= - m_2 \dot \ell(t). 
\end{equation}
In turn, $y^n$ is continuously differentiable and satisfies
\begin{equation} \label{eq:regularised}
    \forall t \geq 0, \: M \dot y^n(t) + \varphi_{1,n_1}'(y^n(t)) + \varphi_{2,n_2}'(y^n(t) + \dot \ell(t)) = - m_2 \ddot \ell(t). 
\end{equation}
Thus, the regularised problem boils down to an ordinary 
differential equation with Lipschitz nonlinear terms depending on $n_1$ and $n_2$. For each $n = (n_1,n_2)$, applying Picard's theorem yields a unique solution $y^n$. 

\subsection{Cauchy convergence of the regularised solutions}

For all $t \in \R_+$, for all $n \in {N^*}^2$ we denote $k_1^n(t) := \int_0^t \varphi_{1,n_1}'(y^n(s)) \textup{d} s$ and $k_2^n(t) := \int_0^t \varphi_{2,n_2}'(y^n(s) + \dot \ell(s)) \textup{d} s + c_2$, with $c_2 = - M y_0 - m_2 \dot \ell(0)$ for compatibility with $(iii)$ in \eqref{myproblem}.
Hereafter, we simply refer to two-indexed sequences as sequences. We say the multi-index $n$ tends to infinity when $\min(n)$ tends to infinity, and that a sequence indexed by $n$ converges if it tends to a limit when $\min n \to \infty$. 

\begin{prop} \label{prop:smooth}
For all $t \geq 0$,
\begin{enumerate}
    \item[$(i)$] The sequence $(y^n)_n$ converges uniformly over $[0,t]$ to a continuous function $y$.
    \item[$(ii)$]  Up to extraction, $k_1^n$ and $k_2^n$ converge pointwise over $[0,t]$ to continuous functions $k_1$, $k_2$. 
    \item[$(iii)$] $y, k_1, k_2$ satisfy the first three lines of system~\eqref{myproblem} over $[0,t]$.
\end{enumerate}
\end{prop}

\begin{proof}[Proof of $(i)$]
Fix $t \geq 0$. Let us start by proving that the set of solutions $\left\{y^{n_1,n_2} \; | \; n_1,n_2\in \mathbb{N}\right\}$ satisfies the Cauchy condition: for all $\varepsilon > 0$ there exists $N > 0$ such that for all $n_1, n_2, r_1, r_2 \in \N^*$, if $\min ( n_1, n_2, r_1, r_2 ) > N$, then $\|y_1^{n_1, n_2} - y_1^{r_1, r_2}\|_{L^\infty([0,t])} < \varepsilon$.
First note that for any function $f$ differentiable over an interval $[0,t]$, we have 
\[
f(t)^2 = f(0)^2+\int_0^t2f(r) \frac{{\rm d} f(r)}{{\rm d}r}\;{\rm d}r
\]
by the fundamental theorem of calculus. Take arbitrary positive naturals  $r_1, r_2, n_1, n_2$, and set $n=(n_1, n_2),\; r = (r_1, r_2)$. Since $y^n (0) = y^r(0)$, injecting equation~\eqref{eq:regularised} into the above gives us
\begin{multline*}
\left( y^n(t) - y^r(t) \right)^2 = - \frac 2 M \int_0^t   \left( y^n(s) - y^r(s) \right) \left[  \varphi_{1,n_1}' (y^n(s)) - \varphi_{1,r_1}' (y^r(s)) \right] \;{\rm d}s \\
    - \frac 2 M \int_0^t   \left( y^n(s) - y^r(s) \right) \left[ \varphi_{1,n_2}' (y^n(s) + \dot \ell (s)) - \varphi_{1,r_2}' (y^r(s) + \dot \ell (s)) \right] \;{\rm d}s,
\end{multline*}
which we can recast into
\begin{multline*}
\left( y^n(t) - y^r(t) \right)^2 = - \frac 2 M \int_0^t   \left( y^n(s) - y^r(s) \right) \left[ \varphi_{1,n_1}' (y^n(s)) - \varphi_{1,r_1}' (y^r(s)) \right] \;{\rm d}s \\
    - \frac 2 M \int_0^t \left( (y^n(s)+\dot \ell (s)) - (y^r(s) + \dot \ell (s)) \right) \left[ \varphi_{1,n_2}' (y^n(s) + \dot \ell (s)) - \varphi_{1,r_2}' (y^r(s) + \dot \ell (s)) \right] \;{\rm d}s.
\end{multline*}
This allows us to apply the following inequality, a consequence of Equation (6.28.b) page 552  in \cite{pardoux2014stochastic}.
\begin{lemma} \label{lem:ysigma}
$\forall y_1, y_2\in \R$, for $i \in \{1, 2\}$ and for all $n_i, r_i \in \N^*$, we have
\[
-(y_1-y_2)(\varphi_{i,{n_i}}'(y_1)-\varphi_{i,{r_i}}'(y_2))\leq (f_i)^2 \left( \frac{1}{n_i}+\frac{1}{r_i} \right).
\]
\end{lemma}
Hence, for all $t \geq 0$, 
\[
\sup_{s\in[0,t]}\left|y^{n}(s)-y^{r}(s)\right|^2 \leq f_1^2 \left(\frac{1}{n_1}+\frac{1}{r_1}\right)t + f_2^2 \left(\frac{1}{n_2}+\frac{1}{r_2}\right) t.
\]
This means that for all $\varepsilon>0$ there exists $N_\varepsilon \geq \frac{2(f_1^2+f_2^2)t}{\varepsilon^2}$ such that for all $r_1, r_2, n_1, n_2 \geq N_\varepsilon$,
\[
\sup_{s\in[0,t]}\left|y^{n}(s) - y^{r}(s)\right|^2\leq\varepsilon^2,
\]
meaning that $\{y^{n}| n_1,n_2\in \mathbb{N}^*\}$ satisfies the Cauchy property. 
It follows that $y^n$ converges uniformly over compact sets to $y \in C^0$ as $n$ tends to $\infty$, concluding the proof of $(i)$.
\end{proof}


\begin{proof}[Proof of (ii)]
    Fix $t \geq 0$. For all $n \in {\N^*}^2$, by definition of $\varphi_{1,n_1}'$ and $\varphi_{2,n_1}'$, the absolute values of the integrands of $k_1^n$ and $k_2^n$ are bounded by $f_1$ and $f_2$ respectively. Hence $(\varphi_{1,n_1}'(y^n))_n$ and $(\varphi_{2,n_1}'(y^n + \dot \ell))_n$ are two sequences bounded in $L^2([0,t])$ by $f_1^2 t$ and $f_2^2 t$ respectively, uniformly in $n$.
    It follows that $(\varphi_{1,n_1}'(y^n))_n$ and  $(\varphi_{2,n_1}'(y^n + \dot \ell))_n$ converge weakly up to extraction to $\delta_1, \delta_2 \in L^2([0,t])$. We may take the same extraction for both sequences by extracting first for one, then for the other.
    Since $\delta_1$ and $\delta_2$ are in $L^2$,
    the functions defined by $k_1(t) := \int_0^t \delta_1 (s) \;{\rm d}s$ and $k_2(t) := \int_0^t \delta_2 (s) \;{\rm d}s$ are continuous.
\end{proof}

\begin{proof}[Proof of $(iii)$]
    For all $n$, equation \eqref{eq:smoothintegrated} is satisfied. 
    For any test function $\phi \in L^2$, up to extraction, 
    \[
\int_0^t \varphi_{1,n_1}'(y^n(s)) \phi(s) \textup{d} s 
+ \int_0^t \varphi_{2,n_2}'(y^n(s) + \dot \ell(s)) \phi(s) \textup{d} s
\xrightarrow[n \to \infty]{} \int_0^t \delta_1(s) \phi(s) \;{\rm d}s + \int_0^t \delta_2(s) \phi(s) \;{\rm d}s.
    \]
So $\phi \equiv 1$ yields, up to extraction, $\int_0^t \varphi_{1,n_1}'(y^n(s)) \textup{d} s + \int_0^t \varphi_{2,n_2}'(y^n(s) + \dot \ell(s)) \textup{d} s \to k_1(t) + k_2(t)$. Since $y^n$ converges uniformly in $C^0$, it follows that the equality  \eqref{eq:smoothintegrated} passes to the limit in $n$ up to extraction, and we recover, for all $t \geq 0$, $M y (t) + k_1(t) + k_2(t) = -m_2 \dot \ell(t)$, concluding the proof of Proposition~\ref{prop:smooth}.

\end{proof}

\subsection{Identification of the limit as a system of variational inequalities}

\begin{prop} \label{prop:vi}
Let $t \geq 0$. The functions $y, k_1, k_2$ from Proposition~\ref{prop:smooth} satisfy
the last two lines of system~\eqref{myproblem} over $[0,t]$.
\end{prop}

\begin{proof}
    Let $n \in \left(\N^*\right)^2$. We first prove that the fourth line of system~\eqref{myproblem} is satisfied by working on $\varphi_1$.
    It is known, see 4. page 550  in \cite{pardoux2014stochastic}, that for all $\xi \in \R$, for $ J_{1,n_1} (\xi) \triangleq \xi - \frac{1}{n_1}\varphi_{1,n_1}'(\xi)$.
\begin{equation} \label{eq:PR13di}
    \varphi_{1,n_1}'(\xi) \in \partial \varphi_1(J_{1,n_1}(\xi)).    
\end{equation} 
Writing out the inclusion~\eqref{eq:PR13di} at $\xi = y^n(r)$ for some $r \in [0,t]$, we recover, for any $v \in \R$,
\[
\varphi_{1,n_1}(y^n(r))[v - J_{1,n_1}(y^n(r))]+ \varphi_1(J_{1,n_1} (y^n(r))) \leq \varphi_1(v).
\]
Furthermore, if $v$ is a continuous function on $[s,t]$, we can integrate the equation above over $[s,t]$: 
\begin{equation} \label{eq:int23}
\int_s^t \varphi_{1,n_1}(y^n(r))[v(r) - J_{1,n_1}(y^n(r))] \;{\rm d}r + \int_s^t \varphi_1(J_{1,n_1} (y^n(r))) \;{\rm d}r \leq  (t-s) \int_s^t \varphi_1(v(r)) \textup{d} r.  \end{equation}

To pass to the limit as $n \to \infty$, 
we rewrite the first term on the left-hand side
\[
\int_s^t \varphi_{1,n_1}(y^n(r))[v - J_{1,n_1}(y^n(r))] \;{\rm d}r 
= \int_s^t \varphi_{1,n_1}(y^n(r))[v - y(r)] \;{\rm d}r 
+ \int_s^t \varphi_{1,n_1}(y^n(r))[y(r) - J_{1,n_1}(y^n(r))] \;{\rm d}r.
\]
Since $y^n$ converges uniformly to $y \in C^0$, $\varphi_{1,n_1}$ is bounded over compact sets, and $J_{1,n_1}(y^n)$ converges uniformly over compact sets to $y$, it follows that the second right-hand side term tends to $0$, and the first converges to
\[
\int_s^t [v(r) - y(r)] \;{\rm d}k_1(r),
\]
by definition of $k_1$.
Since $\varphi_1$ is continuous and $J_{1,n_1}(y^n)$ converges uniformly over compact sets to $y$, we have, for all $r \in [s,t]$,
\[
\varphi_1 (y(r)) \leq \lim_{n\to\infty} \inf \varphi_1 (J_{1,n_1}(y^n(r))).
\]
Integrating over $[s,t]$ and applying Fatou's lemma -- since $\varphi \geq 0$ -- yields
\[
\int_s^t \varphi_1 (y(r)) \;{\rm d} r \leq \int_s^t \lim_{n\to\infty} \inf \varphi_1 (J_{1,n_1}(y^n(r))) \;{\rm d}r \leq \lim_{n\to\infty} \inf \int_s^t \varphi_1 (J_{1,n_1}(y^n(r))) \;{\rm d}r.
\]
Taking the limit inferior of equation~\eqref{eq:int23} and considering the bounds above 
gives us the fourth line of system~\eqref{myproblem}, 
and an analogous reasoning on $\varphi_2$ allows us to recover the fifth line.


\end{proof}

\subsection{Uniqueness of the solution}

Assume two triples of functions $(y,k_1,k_2)$ and $(\tilde y, \tilde{k_1},\tilde{k_2})$ satisfy~\eqref{myproblem}. Then all six functions are continuous, and also $\mathcal{BV}_{loc}$ thanks to the third line of the system, hence differentiable almost everywhere.
Let $0 = s \leq t \in \R$. We can take the variational test function $\mathfrak z$ equal to $\tilde y$ in~\eqref{myproblem}.$(iv)$ 
for $y$ and get
\[
 \int_0^t [\tilde y(r) - y(r)] \textup{d} k_1(r) + \int_0^t \varphi_1(y(r)) \textup{d} r \leq \int_0^t\varphi_1 (\tilde y(r)) \textup{d}r.
\]
Conversely, injecting $\mathfrak z = y$ in in~\eqref{myproblem}.$(iv)$ for $\tilde y$,
\[
 \int_0^t [y(r) - \tilde y(r)] \textup{d} \tilde{k_1}(r) + \int_0^t \varphi_1(\tilde y(r)) \textup{d} r \leq \int_0^t \varphi_1 (y(r)) \textup{d}r.
\]
By addition,
\begin{equation} \label{eq:myproblem4comparison}
    \int_0^t (\tilde y - y)(r)(k_1'-\tilde{k_1}')(r) \;{\rm d}r \leq 0,
\end{equation}
where the derivatives of $k_1$ and $\tilde{k_1}$ are defined weakly. 
Taking $\mathfrak z$ equal to $\tilde y + \dot \ell$ or $y + \dot \ell$ in the fifth lines of~\eqref{myproblem} applied to $y$ and $\tilde y$ respectively then summing like above, we get
\begin{equation} \label{eq:myproblem5comparison}
    \int_0^t (\tilde y - y)(r)(k_2'-\tilde{k_2}')(r) \;{\rm d}r \leq 0.
\end{equation}
Expressing $\tilde y - y$ from~\eqref{myproblem}.$(iii)$ gives us $\tilde y - y = \frac 1 M \left[ (k_1 + k_2) - (\tilde{k_1} + \tilde{k_2}) \right]$, which we inject into the sum of~\eqref{eq:myproblem4comparison} and~\eqref{eq:myproblem5comparison}. Since $M > 0$, we get
\[
\int_0^t \left[ (k_1 + k_2) - (\tilde{k_1} + \tilde{k_2}) \right](r) \left[ (k_1 + k_2) - (\tilde{k_1} + \tilde{k_2}) \right]'(r) \;{\rm d}r \leq 0,
\]
hence
$
\left[ k_1(t) + k_2(t) - \tilde{k_1}(t) - \tilde{k_2}(t) \right]^2 \leq 0,
$
since the expression above vanishes at $t = 0$ thanks to~\eqref{myproblem}.$(iii)$.
It follows that $k_1 + k_2 = \tilde{k_1} + \tilde{k_2}$ at every time, so $\tilde y = y$ thanks to~\eqref{myproblem}.$(iii)$. Hence, subtracting line 4 of~\eqref{myproblem} for $\tilde y$ from the same row for $y$ gives us, for any $0 \leq s \leq t$ and any continuous function $\mathfrak z = y + v$,
\[
\int_s^t v(r) (k_1'(r) - \tilde{k_1}'(r)) \;{\rm d}r \leq 0,
\]
equality resulting from taking $-v$ instead of $v$. And since $k_1(0) = \tilde{k_1}(0)$, we have $k_1 = \tilde{k_1}$, hence $k_2 = \tilde{k_2}$.

\section*{Acknowledgement}
L.M. is thankful for support through NSFC Grant No. 12271364 and GRF Grant No. 11302823.


\bibliographystyle{plain} 
\bibliography{biblio} 

\end{document}